\input amstex
\documentstyle{amsppt}
\input bull-ppt
\keyedby{smit/amh}

\let\sl\it
\define\Que{{\Bbb Q}}
\define\Zee{{\Bbb Z}}
\define\ff{{\Bbb F}}
\define\gp{{\goth p}}
\define\al{\alpha}
\let\goth\germ


\topmatter
\cvol{31}
\cvolyear{1994}
\cmonth{October}
\cyear{1994}
\cvolno{2}
\cpgs{213-215}
\title Zeta functions do not determine class 
numbers\endtitle
\shorttitle{ 
Zeta Functions and Class Numbers}
\author Bart de Smit and Robert Perlis\endauthor
\address Vakgroep Wiskunde, Econometrisch Instituut,
Erasmus Universiteit Rotterdam,  Post\-bus 1738,
3000 DR Rotterdam, Netherlands\endaddress 
\ml 
dsmit\@wis.few.eur.nl\endml 
\address Department of Mathematics, Louisiana State 
University, Baton  
Rouge, Louisiana 70803\endaddress 
\ml 
perlis\@marais.math.lsu.edu\endml 
\date December 1, 1993\enddate
\abstract 
We show that two number fields with the same zeta 
function, and even
with isomorphic adele rings, do not necessarily have the 
same class number.
\endabstract 

\endtopmatter

\document

\heading Introduction\endheading 
The Dedekind zeta function of a number field encodes many  
of its arithmetic invariants, including its degree, 
discriminant,  
number of roots of unity, number of real and complex 
embeddings,
and for every rational
prime the list of residue degrees of the extension primes.
The zeta function also determines the product of the class 
number and
the regulator.  See \cite{3, 4} for proofs.  

It has been an open problem whether or not
the zeta function of a number field
determines its class number.

\proclaim{Theorem} There exist two number fields with the 
same zeta
function and different class numbers.
\endproclaim

Two number fields with the same zeta function are said to 
be {\sl
arithmetically equivalent.}
Let $H$ and $H'$ be subgroups of the Galois group $G$ of a 
Galois
extension $N$ of $\Que$, and denote by $1_H^G$ the 
character of $G$
induced by the trivial character of $H$.
Then the invariant fields $K=N^H$ and $K'=N^{H'}$ are 
arithmetically 
equivalent if and only if $1_H^G=1_{H'}^G$.
If $1_H^G=1_{H'}^G$ and $p$ is a prime not dividing the 
order of $G$,
then the $p$-parts of the class groups of $K$ and $K'$ are 
isomorphic
\cite{4}.

Explicit examples of non-isomorphic arithmetically 
equivalent fields
are the fields $K=\Que({\root 8 \of a})$ and 
$K'=\Que({\root 8\of {16a}})$,
where $a$ is an integer for which none of the numbers $a, 
-a, 2a, -2a$
is a square (see \cite{4}).
The Galois closure of $K$ and $K'$ has degree $32$,
so the class number quotient $h/h'$ is a power of $2$.

\heading Numerical evidence\endheading 

According to our computations with
the system for computational number theory Pari/GP 
(version 1.38,
1993) by Henri Cohen \cite{1, Appendix A}, the values 
$a= -15$, $-31$, $-33$, $-63$, $65$, $66$ give pairs of 
fields
$K$, $K'$ with  $h/h'=2$.
For $a=-65$, $-66$ we get $h'/h=2$.
Wieb Bosma has verified this using the Magma/KANT system 
in Sydney.

Of course this does not qualify as a proof of the theorem.
Perhaps rodents in the bowels of the computer center are 
chewing on wires \pagebreak
and\ altering data.
More disconcerting is the fact that the correctness of the 
algorithms
we used is based on unproven hypotheses; see \cite{1}.

\heading Sketch of proof\endheading 

Take $a=-15$, which gives arithmetically equivalent fields
$K=\Que(\al)$ with $\al^8=-15$ and $K'=\Que(\beta)$ with 
$\beta^8=-240$.

We define the groups of units
$U_0^{\phantom{i}}=\langle -1, {\al+1 \over\al-1}, {\al^2+
\al+2 \over \al+1},
{\al^2-\al+2\over -\al+1}\rangle$
and
$U_0'=\langle -1,$
${\beta^6 + 2\beta^4 - 4\beta^2 -56 \over 16},$
${\beta^7-2\beta^6+2\beta^5-4\beta^3+8\beta^2-8\beta+
64\over 64},$
${\beta^7-2\beta^5+4\beta^4-4\beta^3-32\beta^2+
8\beta-16\over 64}\rangle.$
The fractions were listed by Pari as fundamental units.
We only use that they are units,
which can be checked with a straightforward computation.

The regulators of $U_0^{\phantom{i}}$ and $U_0'$ are 
computed
to be $R_0^{\phantom{i}}\approx 66.316$ and $R_0 '\approx 
132.633$.
Universal regulator bounds from \cite{2} show that the
regulators $R$ and $R'$ of the full unit groups $U$ and $U'$
are at least $0.296$.
It follows that $i=[U:U_0^{\phantom{i}}]$ and 
$i'=[U':U_0']$ are at most
${133\over 0.296} < 500$.
By \cite{3} we have $hR = h'R'$, so 
$R_0'/R_0^{\phantom{i}}= i'/i \cdot h/h'$.
If $h/h'\in\Zee$, then it follows that the denominator of 
the rational number
$R_0'/R_0^{\phantom{i}}$ is at most $500$.
If $h/h'\not\in\Zee$, then $h'/h\in\Zee$, so the denominator
divides $ih'/h = i'R_0^{\phantom{i}}/R_0'<500$.
Our approximations of $R_0^{\phantom{i}}$ and $R_0'$ show 
that
$|2-R_0'/R_0^{\phantom{i}}|<10^{-3}$, and we deduce
that $R_0'/R_0^{\phantom{i}}=2$.

For a unit $u\in U'$ and a prime $\gp$ of $K'$, define
$q(u,\gp)\in\ff_2$ by letting $(-1)^{q(u,\gp)}$ be
the quadratic residue of $u$ modulo $\gp$.
Letting $u$ range over the four generators of $U_0'$
and $\gp$ over the four prime ideals
$(3,\beta)$, $(19,\beta-8)$, $(23,\beta-9)$, 
$(47,\beta-16)$,
one can check that the $(4\times 4)$-matrix 
$(q(u,\gp))_{u,\gp}$
over $\ff_2$ is non-singular.
It follows that $U_0'\cap {U'}^2={U_0'}^2$ and that $i'$ 
is odd.

We now know that $h/h'=R_0'/R_0^{\phantom{i}}\cdot 
i/i'=2i/i'$ and that
$i'$ is odd.
Counting factors $2$, we see that the $2$-power $h/h'$ is 
at least $2$.
This proves the theorem.

With four suitably chosen
primes of $K$ one can also show that $i$ is odd and that 
$h/h'=2$.

\heading Adele rings\endheading 

For each integer $a$ and each odd prime $p$ the 
$\Que_p$-algebras
$\Que_{\, p\, }[X]/(X^8-a)$ and $\Que_{\, p\,  
}[X]/(X^8-16a)$
are isomorphic.
If $a \equiv -1$ modulo $32$, then this also holds for 
$p=2$,
so the fields $K=\Que({\root 8 \of a})$ and 
$K'=\Que({\root 8 \of {16a}})$
have isomorphic adele rings.
For $a=-33$ computations with Pari indicate that the class 
numbers
are different, and this can be verified by imitating the 
argument above.
It follows that the class number of a number field is not 
determined by
the isomorphism class of its adele ring.

\heading Acknowledgment\endheading
The authors are grateful to Wieb Bosma for helpful
discussions and computational assistance.

\enddocument